\documentclass[10pt,a4paper]{amsart}
\usepackage{amsmath,amssymb,mathrsfs,amsfonts,mathtools}
\usepackage{graphicx,geometry}
\usepackage{dynkin-diagrams}
\usetikzlibrary{decorations.pathreplacing,fit,positioning}
\usepackage{setspace}

\usepackage{physics}                    
\usepackage{titleps}

\usepackage{tikz}
\usepackage{tikz-cd}
\usepackage{amssymb}
\usepackage{verbatim}
\usepackage{textcomp}

\usepackage{float,comment,xcolor}
\usepackage{enumitem}
\usepackage{tikz-cd}
\usepackage{fancyhdr}
\newcommand{\shorttitle}{}

\graphicspath{{plot/}}
\newtheorem{theorem}{Theorem}[section]
\newtheorem{lemma}[theorem]{Lemma}

\newtheorem{corollary}[theorem]{Corollary}

\theoremstyle{definition}
\newtheorem{definition}[theorem]{Definition}

\theoremstyle{remark}

\numberwithin{equation}{section}
\newcommand{\abbrevtitle}{}
\def\nd{\noindent}

\title[Gauss-Bonnet formula for metrics with logarithmic singularities]{Gauss-Bonnet formula for metrics with logarithmic singularities on compact Riemann surfaces}

\author{Yuanjiu Lyu$^\dagger$}
\address{School of Mathematical  Sciences, University of Science and Technology of China, Hefei 230026 China}
\email{lyj999@mail.ustc.edu.cn}

\author{Bin Xu}
\address{CAS Wu Wen-tsun Key Laboratory of Mathematics and School of Mathematical  Sciences, University of Science and Technology of China, Hefei 230026 China}
\email{bxu@ustc.edu.cn}

\thanks{B.X. is supported in part by the Project of Stable Support for Youth Team in Basic Research Field, CAS (Grant No. YSBR-001) and NSFC (Grant Nos. 12271495 and 12071449). }

\thanks{$^\dagger$Y.L. is the corresponding author.}

\fancypagestyle{plain}{
  \fancyhead[LE,RO]{\shorttitle}
  \fancyfoot[C]{\thepage}
}

\begin{document}
\pagestyle{plain}

\maketitle

\begin{abstract}
We prove a generalization of the classical Gauss-Bonnet formula for metrics with logarithmic singularities on compact Riemann surfaces, under the condition that the Gaussian curvature is Lebesgue integrable with respect to the metric's area form. As an application, we establish the Gauss-Bonnet formula for special K\"ahler metrics when their associated cubic differentials are meromorphic on compact Riemann surfaces.    

\vskip 4.5mm

\nd \begin{tabular}{@{}l@{ }p{10.1cm}} {\bf Keywords } &

Gauss-Bonnet formula; Logarithmic singularity; Riemann surface; special K\"ahler metric.
\end{tabular}

\nd {\bf 2020 MR Subject Classification } 
53C21, 14H15

\end{abstract}

\section{Introduction}

S.S. Chern \cite{Chern_intrinsic-proof-GaussBonnet1944} provided an intrinsic proof of the Gauss-Bonnet formula. Specifically, for any compact orientable smooth Riemannian manifold of even dimension, the Euler class is expressible in terms of the curvature of a smooth Riemannian metric.  
For Riemannian metrics with edge-cone singularities along an embedded 2-manifold, M. Atiyah and C. Lebrun \cite{Atiyah-Lebrun2013} established modified versions of the Gauss-Bonnet and signature theorems. Additionally, C. T. McMullen \cite{Mcmullen_GaussBonnet-cone-manifolds2017} proved a Gauss-Bonnet analogue for cone manifolds.  
On a compact Riemann surface $X$, M. Troyanov \cite{Troyanov_compact-surface-conical-singularities1991} generalized the classical Gauss-Bonnet formula to Riemannian metrics with conical singularities, assuming the Gaussian curvature extends to a H\"older continuous function on $X$.  
Subsequently, Fang-Xu-Yang \cite{Gauss-Bonnet_FXY2024} relaxed the H\"older condition to $L^1$-integrability of the Gaussian curvature $K$ with respect to the metric's area form.

In this paper, we generalize Fang-Xu-Yang's results \cite{Gauss-Bonnet_FXY2024} to a broader class of singular metrics (Theorem \ref{gauss-bonnet-general}).  
As an application, we investigate special Kähler metrics arising in mathematical physics (see, e.g., Freed \cite{spK-manifolds_Freed1900}).  
A. Haydys and the second author \cite{Haydys-Xu_sks2020} derived necessary conditions for the existence of special K\"ahler structures with meromorphic cubic forms on punctured Riemann surfaces, establishing sufficiency for the sphere and tori. We prove that the Gaussian curvature of any special K\"ahler metric with meromorphic cubic form is $L^1$-integrable, which implies the Gauss-Bonnet formula for such metrics via Corollary \ref{Cor_GaussBonnet_sks}.

We now introduce the framework for our main results. Let $X$ be a compact Riemann surface and $D = \{p_1,\dots,p_n\}\subset X$ a finite set. For each $j = 1,\dots,n$, choose a chart $(V_j,z_j)$ centered at $p_j$ and a Euclidean disk $U_j \subset V_j$ centered at $p_j$ such that $U_i \cap U_j = \emptyset$ for $i \neq j$.  
Consider a smooth K\"ahler metric $g$ on $X \setminus D$ with the following local property: in local coordinates $z_j$ on each $U_j$, $g = e^{2u} |\mathrm{d} z_j|^2$ for some smooth function $u$ on $U_j\setminus\{p_j\}$ satisfying
\begin{equation}
\label{local_u}
    u = \alpha_j \log r + \sum_{\ell=1}^{k_j} \beta_{\ell j} \log^{(\ell)} r + v_j, \quad r = |z_j|; \,\, \alpha_j, \beta_{\ell j} \in \mathbb{R}
\end{equation}
where $v_j$ is bounded and smooth on $U_j \setminus \{p_j\}$, and $\log^{(\ell)} r$ denotes the $\ell$-fold iterated logarithm (i.e., $\log^{(\ell)} r := \underbrace{\log \cdots \log}_{\ell} |\log r|$). When $\beta_{\ell j} = 0$ for all $\ell,j$ and each $v_j$ extends continuously to $p_j$, $g$ is a conical metric in the sense of Fang-Xu-Yang \cite{Gauss-Bonnet_FXY2024}.

\begin{definition}
    The \emph{order of $g$ at $p_j$} is well-defined as
\begin{equation*}
    \mathrm{ord}(g, p_j) := \alpha_j.
\end{equation*}
Denote by $K$ the Gaussian curvature of $g$ and $\mathrm{d}A$ the associated area form on $X \setminus D$. 
\end{definition}

Our main results are stated as follows:

\begin{theorem}
\label{gauss-bonnet-general}
    With the preceding notation, if $K$ is $L^1$-integrable with respect to $\mathrm{d}A$ (i.e., $\int_{X \setminus D} |K|  \mathrm{d}A < \infty$), then
    \begin{equation}
    \label{eq:gauss-bonnet-general}
        \frac{1}{2\pi} \int_{X \setminus D} K  \mathrm{d}A = \chi(X) + \sum_{j=1}^n \mathrm{ord}(g, p_j),
    \end{equation}
    where $\chi(X) = 2 - 2g(X)$ is the Euler characteristic of $X$ and $g(X)$ its genus.
\end{theorem}

\begin{corollary}
\label{Cor_GaussBonnet_sks}
    Let $g$ be a special K\"ahler metric on $X \setminus D$ whose associated cubic form is meromorphic on $X$. Then $K$ is $L^1$-integrable with respect to $\mathrm{d}A$. Consequently, Theorem \ref{gauss-bonnet-general} yields
    \begin{equation}
    \label{eq:gaussbonnet-SKS}
        \frac{1}{2\pi} \int_{X \setminus D} K  \mathrm{d}A = 2 - 2g(X) + \sum_{j=1}^n \mathrm{ord}(g, p_j).
    \end{equation}
\end{corollary}

\section{Proofs of the main results}

In this subsection, we give the proofs of Theorem \ref{gauss-bonnet-general} and Corollary \ref{Cor_GaussBonnet_sks}. In the following, we denote $\Delta = \Delta_{\mathbb{R}^2}$ the standard Laplacian operator on $\mathbb{R}^2$, in the polar coordinate $(r,\theta)$, it is given by 
\begin{equation*}
    \Delta f(r) = \frac{1}{r}\frac{\partial}{\partial r} \left(r\frac{\partial f}{\partial r}\right) = \frac{\partial^2 f}{\partial r^2} + \frac{1}{r} \frac{\partial f}{\partial r}
\end{equation*}
we can also verify that $\Delta f \cdot \mathrm{vol}_{\mathbb{R}^2}= \dd \star \dd f $ where $\star$ is the standard Hodge star operator on $\mathbb{R}^2$.

\subsection{Proof of Theorem \ref{gauss-bonnet-general}} 
We need the following two technical lemmas for preparation,
\begin{lemma}
\label{lem_laplacian-log-k-r}
Assume $r \in(0,1)$. There hold
    \begin{align*}
        &\Delta \log r = 0 \quad \text{ and }\\
        &\Delta \log^{(k)}r = 
        \frac{P_n(\log r, \cdots,\log^{(k-1)}r)}
        {r^2 (\log r)^2 \cdots(\log^{(k-1)}r)^2} \quad  \text{ for } k \geq1. 
        \end{align*}
     where $P_n = -1$ for $k = 1$ and $P_n = - \prod_{j = 1}^{k-1}\log^{(j)}r 
        - \sum_{l = 1}^{k-1}  
        \prod_{j = l+1}^{k-1}\log^{(j)}r$ for $k\geq2$. Let $B_\epsilon(0)$ be a ball of radius $\epsilon<1$ centered at $0\in\mathbb{R}^2$, then for $k\geq1$,
     \begin{equation*}
         \int_{B_\epsilon(0)\setminus0} \left|\Delta\log^{(k)}r\right| <\infty.
     \end{equation*}
\end{lemma}
\begin{proof}
    Denote $\log^{(0)}r := |\log r| = -\log r$. Then for $k\geq1$
    \begin{align*}
        \frac{\partial \log^{(k)}r}{\partial r} 
        &= \frac{\partial }{\partial r}\left(\log \left(\log^{(k-1)}r\right)\right)\\
        &=\frac{1}{\log^{(k-1)}r} \cdot
        \frac{\partial \log^{(k-1)}r}{\partial r}\\
        &=
        -\frac{1}{r\cdot\prod_{j = 0}^{k-1}\log^{(j)}r}
    \end{align*}
    and then
    \begin{align*}
        \text{ for } k = 1, \quad 
        \frac{\partial}{\partial r}\left(\frac{\partial \log^{(1)}r}{\partial r}\right)
        =
        - \frac{\log r + 1}{r^2(\log r)^2} 
    \end{align*}
    \begin{align*}
    \text{ for } k \geq 2, \quad
        \frac{\partial}{\partial r}\left(\frac{\partial \log^{(k)}r}{\partial r}\right)
        &=
        -\frac{\partial}{\partial r}\left(
        \frac{1}{r\cdot\prod_{j = 0}^{k-1}\log^{(j)}r}
        \right)\\
        &=
        \frac{1}{\left(r\cdot\prod_{j = 0}^{k-1}\log^{(j)}r\right)^2} \cdot
        \frac{\partial}{\partial r}\left(
        r\cdot\prod_{j = 0}^{k-1}\log^{(j)}r
        \right)\\
        &=
        \frac{1}{\left(r\cdot\prod_{j = 0}^{k-1}\log^{(j)}r\right)^2} \cdot
        \left(
        \begin{aligned}
            &\frac{\partial}{\partial r}(r)\cdot\prod_{j = 0}^{k-1}\log^{(j)}r\\
            &+ r\cdot\left(\frac{\partial}{\partial r}\log^{(0)}r\right)\cdot\prod_{j=1}^{k-1}\log^{(j)}r\\
            &+ r\log^{(0)} r \cdot \frac{\partial}{\partial r}\left(\prod_{j=1}^{k-1}\log^{(j)}r\right)
        \end{aligned}
        \right)\\
        &=
        \frac{1}{\left(r\cdot\prod_{j = 0}^{k-1}\log^{(j)}r\right)^2} \cdot
        \left(
        \prod_{j = 0}^{k-1}\log^{(j)}r 
        - \prod_{j = 1}^{k-1}\log^{(j)}r 
        + r\log^{(0)} r
        \sum_{l = 1}^{k-1} \left( \frac{\partial \log^{(l)}r}{\partial r} 
        \prod_{\substack{j = 1\\j\neq l}}^{k-1}\log^{(j)}r\right)
        \right)\\
        &=
        \frac{1}{\left(r\cdot\prod_{j = 0}^{k-1}\log^{(j)}r\right)^2} \cdot
        \left(
        \prod_{j = 0}^{k-1}\log^{(j)}r 
        - \prod_{j = 1}^{k-1}\log^{(j)}r 
        - \sum_{l = 1}^{k-1} \left(
        \frac{1}{\prod_{j = 1}^{l-1} \log^{(j)}r}
        \cdot 
        \prod_{\substack{j = 1\\j\neq l}}^{k-1}\log^{(j)}r\right)
        \right)\\
        &=
         \frac{1}{\left(r\cdot\prod_{j = 0}^{k-1}\log^{(j)}r\right)^2} \cdot
        \left(
        \prod_{j = 0}^{k-1}\log^{(j)}r 
        - \prod_{j = 1}^{k-1}\log^{(j)}r 
        - \sum_{l = 1}^{k-1}  
        \prod_{j = l+1}^{k-1}\log^{(j)}r
        \right)
    \end{align*}
    Thus we see
    \begin{align*}
        \Delta \log^{(1)} r 
        = \frac{\partial^2 \log^{(1)}r}{\partial r^2} + \frac{1}{r}\frac{\partial \log^{(1)}r}{\partial r}
        =
        \frac{-1}{r^2(\log r)^2}
    \end{align*}
    and for $k\geq2$,
    \begin{align*}
        \Delta \log^{(k)} r 
        &= \frac{\partial^2 \log^{(k)}r}{\partial r^2} + \frac{1}{r}\frac{\partial \log^{(k)}r}{\partial r}\\
        &=
        \frac{1}{\left(r\cdot\prod_{j = 0}^{k-1}\log^{(j)}r\right)^2} \cdot
        \left(
        \cdot
        \prod_{j = 0}^{k-1}\log^{(j)}r 
        - \prod_{j = 1}^{k-1}\log^{(j)}r 
        - \sum_{l = 1}^{k-1}  
        \prod_{j = l+1}^{k-1}\log^{(j)}r
        \right)
        -
        \frac{1}{r^2\cdot\prod_{j = 0}^{k-1}\log^{(j)}r}\\
        &=
        \frac{1}{\left(r\cdot\prod_{j = 0}^{k-1}\log^{(j)}r\right)^2} \cdot
        \left(
        - \prod_{j = 1}^{k-1}\log^{(j)}r 
        - \sum_{l = 1}^{k-1}  
        \prod_{j = l+1}^{k-1}\log^{(j)}r
        \right)
    \end{align*}
    Then 
    \begin{align*}
        \int_{B_\epsilon(0)\setminus0}\left|\Delta\log^{(1)}r\right|
        =
        \int_{B_\epsilon(0)\setminus0} \frac{1}{r^2(\log r)^2}
        =
        2\pi \int_0^\epsilon \frac{1}{r(\log r)^2} \dd r
        =
        -\frac{1}{\log \epsilon}
    \end{align*}
    and for $k\geq2$,
    \begin{align*}
        \int_{B_\epsilon(0)\setminus0}\left|\Delta\log^{(k)}r\right|
        &=
        2\pi \cdot \int_{0}^\epsilon\left|\Delta\log^{(k)}r\right| r \dd r\\
        &=
        2\pi
        \int_0^\epsilon
        \frac{1}{\left(r\cdot\prod_{j = 0}^{k-1}\log^{(j)}r\right)^2} \cdot
        \left(
        - \prod_{j = 1}^{k-1}\log^{(j)}r 
        - \sum_{l = 1}^{k-1}  
        \prod_{j = l+1}^{k-1}\log^{(j)}r
        \right)r\dd r\\
        &=
        2\pi\int^{\log\epsilon}_{-\infty} 
        \frac{1}{t^2} 
        \cdot
        \left(\frac{1}{\log |t| \cdots \log^{(k-1)}|t|}\right)^2
        \cdot
        \left(
        - \prod_{j = 1}^{k-1}\log^{(j)}|t| 
        - \sum_{l = 1}^{k-1}  
        \prod_{j = l+1}^{k-1}\log^{(j)}|t|
        \right)
        \dd t < \infty
    \end{align*}
    where we set $t = \log r$, and $\log^{(k)}|t| : = \underbrace{\log\cdots\log}_k|t|$ and assume $\epsilon$ is small enough.
\end{proof}

\begin{lemma}
\label{lemma1}
    Let \( f \) be a bounded smooth function defined on \( U_j \setminus \{p_j\} \), more precisely, we assume \( f \) extends smoothly to slightly larger open sets $U_j^k\setminus\{p_j\}$ where $U_j^k := \{z\in\mathbb{C}\mid d(z,U_j)\leq \frac{1}{k}\}$ for $k\in\mathbb{N}$. If \( f \) satisfies
    \begin{equation*}
        \frac{\mathrm{i}}{2}\int_{U_j \setminus \{p_j\}} |\Delta f| \, \mathrm{d} z_j\wedge\dd \Bar{z}_j < \infty,
    \end{equation*}
    then the following holds:
    \begin{enumerate}
        \item Let $\Delta_1 f$ denote the distributional Laplacian of $f$ by viewing $f$ as a distribution in $U_j$. Since $\Delta f$ is smooth on $U_j^k \setminus \{p_j\}$ and satisfies $\int_{U_j \setminus \{p_j\}} |\Delta f| < \infty$, it defines a distribution $\Delta_2 f$ in $U_j$. These two distributions coincide:
        \begin{equation*}
            \Delta_1 f = \Delta_2 f \quad \text{ in } \,U_j.
        \end{equation*}
        \item The Dirichlet energy of \( f \) is finite:
        \begin{equation*}
            \frac{\mathrm{i}}{2}\int_{U_j \setminus \{p_j\}} \|\nabla f\|^2 \mathrm{d} z_j\wedge\dd \Bar{z}_j < \infty,
        \end{equation*}
        where \( z_j = x + \mathrm{i}y \) (we would abbreviate $z_j$ to $z$), \( \nabla f = (\partial_x f, \partial_y f) \), and \( \|\nabla f\|^2 = (\partial_x f)^2 + (\partial_y f)^2 \).

        \item The flux of \( f \) vanishes asymptotically at \( p_j \):
        \begin{equation}
            \liminf_{\epsilon \to 0} \left| \int_{\partial U_j(\epsilon)} \star \mathrm{d} f \right| = 0.
        \end{equation}
    \end{enumerate}
\end{lemma}
\begin{proof}
    (1): The distribution $\Delta_1 f - \Delta_2 f$ has support contained in $\{p_j\}$. By the structure theorem for distributions with point support \cite[pp. 46-47, Theorem 2.3.4]{Hormander_PDE-I-1983}, there exists $q \geq 0$ such that
    \begin{equation*}
        \Delta_1 f - \Delta_2 f = \sum_{|\alpha|\leq q} a_\alpha \partial^\alpha \delta_{p_j}
    \end{equation*}
    where $\delta_{p_j}$ is the Dirac $\delta$-function supported at $p_j$.
    To show all coefficients $a_\alpha$ vanish, fix an arbitrary multi-index $\alpha$ and select $\phi \in C_0^\infty(U_j)$ satisfying:
    \begin{equation*}
        \partial^\alpha \phi(p_j) \neq 0 \quad \text{and} \quad \partial^\beta \phi(p_j) = 0 \quad \text{for all } \beta \neq \alpha \text{ with } |\beta| \leq q.
    \end{equation*}
    Define the scaled test functions $\phi_k(z) := \phi(k\cdot z)$. 
    First, we observe that $\left(\Delta_2 f\right)(\phi_k) \to 0$ as $k\to \infty$ since:
    \begin{align*}
        |\left(\Delta_2 f\right)(\phi_k)| 
        &= \left|\int_{U_j \setminus \{p_j\}} \Delta f(z) \phi_k(z) \right| \\
        &\leq \|\phi\|_{L^\infty} \int_{\mathrm{supp}(\phi_k) \setminus \{p_j\}} |\Delta f(z)| \to 0 \quad \text{as } k \to \infty.
    \end{align*}
    Next, let $c_0 := \limsup_{p \to p_j} f(p) < \infty$ and $c_1 := \liminf_{p \to p_j} f(p) >-\infty$ (which exist by boundedness of $f$). The scaled functions satisfy:
    \begin{equation*}
        \limsup_{k \to \infty} f(z/k) = c_0 \text{ and } 
        \liminf_{k \to \infty} f(z/k) = c_1
        \quad \text{ on } \mathrm{supp}(\phi)\setminus\{p_j\}.
    \end{equation*}
    Applying Fatou's Lemma and reverse Fatou's lemma yields:
    \begin{align*}
        0 = 
        \int_{U_j \setminus \{p_j\}} 
        \liminf_{k \to \infty} (f(z/k) - c_1) \Delta \phi(z)  &\leq 
        \liminf_{k \to \infty} \int_{U_j}(f(z/k) - c_1) \Delta \phi(z) ,\\
        \limsup_{k \to \infty} \int_{U_j \setminus \{p_j\}} (f(z/k) - c_0) \Delta \phi(z)   
        &\leq \int_{U_j} \limsup_{k \to \infty} (f(z/k) - c_0) \Delta \phi(z)  = 0.
    \end{align*}
    Since $(\Delta_1f)(\phi_k) = \int_{U_j \setminus \{p_j\}} f(\cdot/k) \Delta \phi(\cdot) $, we bound $(\Delta_1f)(\phi_k)$ as
    \begin{align*}
        -\infty
        <
        \int_{U_j}c_1\Delta\phi
        \leq
        \liminf_{k\to \infty} (\Delta_1f)(\phi_k)
        \leq
        \limsup_{k\to \infty} (\Delta_1f)(\phi_k)
        \leq 
        \int_{U_j} c_0\Delta\phi
        <\infty
    \end{align*}
    However, the distributional difference gives:
    \begin{equation*}
        (\Delta_1 f - \Delta_2 f)(\phi_k) = a_\alpha k^{|\alpha|} \partial^\alpha \phi(0)
    \end{equation*}
    For this to remain bounded as $k \to \infty$, we must have $a_\alpha = 0$ since $\partial^\alpha \phi(0) \neq 0$. Thus $\Delta_1 f = \Delta_2 f$ in $U_j$.
    
    (2):
    First, we define $f_k = \chi_k*f$, where $\chi_k\in C^\infty_0(\mathbb{C})$ with the property
    \begin{equation*}
        \chi_k \geq 0\,, \,\,\int_\mathbb{C}\chi = 1 \text{ and } \mathrm{supp}(\chi_k) \subset \left\{|z|<\frac{1}{k}\right\}. 
    \end{equation*}
    Then for any $q\in U_j\setminus\{p_j\}$, we have the pointwise convergence: $f(q) = \lim_{k\to \infty}f_k(q)$ and  $\nabla f(q) = \lim_{k\to \infty} \nabla f_k (q)$. Applying Fatou's lemma yields:
    \begin{equation*}
        \int_{U_j\setminus\{p_j\}}\norm{\nabla f}^2  
        \leq 
        \liminf_{k\to \infty } \int_{U_j}\norm{\nabla f_k}^2 
    \end{equation*}
    
    We now claim that there holds a uniform bound $\int_{U_j} \|\nabla f_k\|^2 |\mathrm{d} z|^2 < C$ for all $k \in \mathbb{N}$:
    Since $f_k$ are smooth on $U_j$, Stokes formula yields
    \begin{align*}
        \int_{U_j} \dd \left( 
        f_k\left(\partial_x(f_k)\dd y - \partial_y(f_k)\dd x\right)
        \right)  
        &= \int_{U_j}\norm{\nabla f_k}^2 \dd x\wedge\dd y 
        + 
        \int_{U_j}f_k \cdot \Delta f_k \dd x\wedge \dd y\\
        &= \int_{\partial U_j} f_k\left(
        \partial_x(f_k)\dd y - \partial_y(f_k)\dd x
        \right)
    \end{align*}
    which implies:
    \begin{equation*}
        \int_{U_j}\norm{\nabla f_k}^2 \dd x\wedge\dd y 
        = 
        -\int_{U_j}f_k \cdot \Delta f_k \dd x\wedge \dd y
        +\int_{\partial U_j} f_k\left(
        \partial_x(f_k)\dd y - \partial_y(f_k)\dd x
        \right)
    \end{equation*}
    The second term converges to $\int_{\partial U_j}f\left(\partial_xf\dd y - \partial_yf\dd x\right)$ hence uniformly bounded as $k\to \infty$. The boundedness of the first term follows from $|f_k| < |\chi_k|\cdot|f| < C$ and the boundedness of $\int_{U_j} |\Delta f_k|$ as:
    \begin{align*}
        \int_{U_j}\left|\Delta f_k \right|
        &= \int_{U_j} \left| \Delta\left(\chi_k*f\right)\right|\\
        &= \int_{U_j} \left|\chi_k*\Delta_1 f\right| \\
        &= \int_{U_j} \left|\chi_k*\Delta_2 f\right| \\
        &= \int_{Uj}\left|\int_{U_j^k}\chi_k(z-w)\cdot \Delta_2 f(w) \dd w\right|\dd z\\
        &\leq
        \int_{\mathbb{R}^2}\chi_k \int_{U^k_j} \Delta_2 f < \infty
    \end{align*}
    where we use (1) in the third equality.
    Hence, 
    \begin{equation*}
        \frac{\mathrm{i}}{2}\int_{U_j\setminus\{p_j\}} \norm{\nabla f}^2 \mathrm{d} z\wedge\dd \Bar{z} < C
    \end{equation*}
    
    (3):
    We have shown that $\norm{\nabla f}\in L^2(U_j\setminus\{p_j\})$. In polar coordinates $z = re^{\mathrm{i}\theta}$, we express the Hodge star operation as
    \begin{align*}
        \star \dd f = \star (f_r\dd r + f_\theta\dd \theta) 
        = r(f_r\dd \theta - f_\theta \dd r)
    \end{align*}
    which yields the boundary integral:
    \begin{equation*}
        \int_{\partial U_j(\epsilon)} \star\dd f = \epsilon\int_{r = \epsilon} \partial_rf \dd \theta
    \end{equation*}
    To prove $\liminf_{\epsilon \to 0} \epsilon \int_{0}^{2\pi} |\partial_r f| \, \mathrm{d} \theta = 0$, we proceed by contradiction. Suppose the contrary holds: $\liminf_{\epsilon\to 0}\epsilon\int_{r = \epsilon}|\partial_rf| \dd \theta = C>0$, that is, there exists $\epsilon_0>0$ such that for any $\rho < \epsilon_0$, there holds
    \begin{equation*}
        \int_{r = \rho} |\partial_rf| \dd \theta > \frac{C}{\rho}.
    \end{equation*}
    Denote $U_j(\rho)\subset U_j$ as the disk of Euclidean radius $\rho$, that is, $U_j(\rho) := \{z\in U_j\mid |z|<\rho\}$. Applying the Cauchy-Schwarz inequality yields:
    \begin{align*}
       \pi\rho^2
       \cdot\int_{U_j(\rho)\setminus\{p_j\}} \norm{\nabla f}^2
       =
       \int_{U_j(\rho)\setminus\{p_j\}}1 \cdot\int_{U_j(\rho)\setminus\{p_j\}} \norm{\nabla f}^2
       \geq 
       \left(\int_{U_j(\rho)\setminus\{p_j\}}\norm{\nabla f}\right)^2.
    \end{align*}
    Noting that $\|\nabla f\| = \|\mathrm{d} f\|_{\text{Eucl}}$ (the Euclidean norm of the 1-form), we estimate:  
    \begin{align*}
        \int_{U_j(\rho)\setminus\{p_j\}}\norm{\nabla f} 
        &=
        \int_{U_j(\rho)\setminus\{p_j\}}\norm{\partial_rf\dd r + \partial_\theta\dd \theta}\\
        &\geq
        \int_{U_j(\rho)\setminus\{p_j\}}\norm{\partial_r f\dd r}\\
        &=
        \int_{U_j(\rho)\setminus\{p_j\}} |\partial_rf|\\
        &=
        \int_{0}^{\rho} \left(\int_{} |\partial_rf|  \dd \theta\right) r\dd r\\
        &\geq
        C\rho
    \end{align*}
    Combining these inequalities yields:
    \begin{align*}
        \int_{U_j(\rho)\setminus\{p_j\}} \norm{\nabla f}^2
        \geq 
        \frac{C^2}{\pi}\quad \text{for any } \,\,\rho<\epsilon_0
    \end{align*}
    However, this contradicts the fact that for any $L^1$ function $F$:
    \begin{equation*}
        \lim_{n\to \infty} \int_{A_n} |F| \to 0 \quad \text{if} \quad \lim_{n\to\infty}\mu(A_n) = 0.
    \end{equation*}
    Thus 
    \begin{align*}
        \liminf_{\epsilon\to0} \left|\int_{\partial U_{j}(\epsilon)} \star \dd f\right|
        \leq 
        \liminf_{\epsilon\to0} \epsilon\int_{r = \epsilon} |\partial_rf| \dd \theta = 0
    \end{align*}
\end{proof}

Now we can prove Theorem \ref{gauss-bonnet-general}:
Fix a smooth K\"ahler metric $g_0$ on $X$ with Gaussian curvature $K_0$ and area form $\mathrm{d}A_0$. Then there exists a smooth function $v$ on $X \setminus D$ such that 
\[
    g = e^{2v} g_0 \quad \text{on } X \setminus D,
\]
where $v$ is determined by the conformal factor relating the area forms ($e^{2v} = \mathrm{d}A / \mathrm{d}A_0$). The curvature forms satisfy the relation
\[
    K  \mathrm{d}A = K_0  \mathrm{d}A_0 - \mathrm{d}(\star \mathrm{d} v),
\]
where $\star$ denotes the Hodge star operator with respect to the complex structure of $X$ (which is metric-independent). 
The classical Gauss-Bonnet theorem yields:
\[
    \frac{1}{2\pi} \int_X K_0  \mathrm{d}A_0 = \chi(X) = 2 - 2g(X),
\]
where $\chi(X)$ is the Euler characteristic of $X$ and $g(X)$ its genus. 
Therefore, to establish (\ref{gauss-bonnet-general}), it suffices to prove
\[
    -\frac{1}{2\pi} \int_{X \setminus D} \mathrm{d}(\star \mathrm{d} v) = \sum_{j=1}^n \mathrm{ord}(g, p_j).
\]

    In the local coordinate system $(U_j, z_j)$, the metrics are expressed as $g = e^{2u} |\mathrm{d} z_j|^2$ and $g_0 = e^{2u_0} |\mathrm{d} z_j|^2$, giving $v = u - u_0$ on $U_j\setminus\{p_j\}$. 
    Denote $U_j(\epsilon)\subset U_j$ as the disk of Euclidean radius $\epsilon<1$, that is, $U_j(\epsilon) := \{z\in U_j\mid |z_j|<\epsilon\}$. Applying Stokes' theorem yields:
\begin{equation*}
    \frac{1}{2\pi} \int_{X \setminus D} \mathrm{d}(\star \mathrm{d} v)
    = \lim_{\epsilon \to 0} 
    \frac{1}{2\pi} \int_{X \setminus \bigcup_j U_j(\epsilon)} \mathrm{d}(\star \mathrm{d} v)
    = -\lim_{\epsilon \to 0} \sum_j
    \frac{1}{2\pi} \int_{\partial U_j(\epsilon)} \star \mathrm{d} v.
\end{equation*}
Since $u_0$ is smooth on $U_j$, we observe as $\epsilon \to 0$:
\begin{equation*}
    \int_{\partial U_j(\epsilon)} \star \mathrm{d} u_0 \to 0.
\end{equation*}
Therefore, it suffices to compute:
\begin{equation*}
    \liminf_{\epsilon \to 0} \frac{1}{2\pi} \int_{\partial U_j(\epsilon)} \star \mathrm{d} u
\end{equation*}
for each $j$. Substituting the expression for $u$:
\begin{align*}
    \frac{1}{2\pi} \int_{\partial U_j(\epsilon)} \star \mathrm{d} u
    &= \frac{1}{2\pi} \int_{\partial U_j(\epsilon)} \star \mathrm{d} \left( 
        \alpha_j \log r + \sum_{\ell=1}^{k_j} \beta_{\ell j} \log^{(\ell)} r + v_j 
    \right) \\
    &= \alpha_j + \frac{1}{2\pi} \int_{\partial U_j(\epsilon)} \star \mathrm{d} v_j + \mathcal{O}\left( \frac{1}{\log \epsilon} \right).
\end{align*}

      The condition $\int_{U_j\setminus\{p_j\}}|K|\dd A < \infty$ implies
        \begin{equation*}
        \int_{U_j\setminus\{p_j\}} |\Delta u| <\infty.
        \end{equation*}
      By calculation of Lemma \ref{lem_laplacian-log-k-r}, 
     \begin{align*}
         \int_{U_j\setminus0}\Delta  \log r = 0 < \infty\quad \text{and}\quad
         \int_{U_j\setminus0} |\Delta \log^{(k)}r| < \infty 
     \end{align*}
     Hence, combining with $\int_{U_j\setminus\{p_j\}}|\Delta u|<\infty$ and the decomposition $u = \alpha_j\log r + \sum_{\ell}\beta_{\ell j}\log^{(\ell)}r + v_j$, there holds
    \begin{equation*}
        \int_{U_j\setminus\{p_j\}} |\Delta v_j| <\infty.
    \end{equation*}
    Then by the Lemma \ref{lemma1} above, as $\epsilon\to 0$, we see
    \begin{equation*}
    \liminf_{\epsilon\to0}\frac{1}{2\pi} \int_{\partial U_j(\epsilon)} \star \dd v_j = \alpha_j = \mathrm{ord}(g,p_j).
    \end{equation*}

\qquad \qquad\qquad\qquad\qquad\qquad\qquad\qquad \qquad \qquad\qquad\qquad\qquad\qquad\qquad\qquad \qquad \qquad \qquad
\qquad \qquad
\boxed{}
\newline

\subsection{Proof of Corollary \ref{Cor_GaussBonnet_sks}} Firstly, there are only two models for local special K\"ahler metrics:
\begin{lemma}(A. Haydys, \cite{isolated_SKS-Haydys2015})
\label{A._local_model_SKS}
    Let $g_{\mathrm{local}}$ be an special K\"ahler metric on the punctured unit disk $B_1^*(0)\subset\mathbb{C}$. We assume the associated cubic form $\Theta$ is meromorphic over $B_1(0)$ with a pole of order $n >-\infty$ at the origin. 
    Denote $z$ as the natural coordinate of $B_1(0)\subset\mathbb{C}$, then
    \begin{align*}
        g_{\mathrm{local}} &= - e^{O(1)}\cdot r^{n+1} \log r |\dd z|^2\\
        &or\\
        g_{\mathrm{local}} &= r^\beta (C + o(1)) |\dd z|^2
    \end{align*}
    where $C >0$ and $\beta < n+1$.
    
\end{lemma}

Therefore, the special K\"ahler metric $g$ with meromorphic condition is contained in the class of singular metrics in subsection 1.1 and the order of $g$ at $p_j$ is given by:
\begin{equation*}
    \mathrm{ord}(g,p_j) =
    \begin{cases}
        \frac{n_j+1}{2} & \text{if } g = - e^{O(1)}\cdot r^{n_j+1} \log r |\dd z_j|^2\\
        \frac{\beta_j}{2} & \text{if } g = r^{\beta_j} (C + o(1)) |\dd z_j|^2
    \end{cases}
\end{equation*}
where $z_j$ is some coordinate on $X$ centered at $p_j$ as before.
Thus, to establish the Gauss-Bonnet formula, we only need to show the $L^1$ integrability of its Gaussian curvature.
    Write $g = e^{-U}|\dd z_j|^2$ and the associated cubic form $\Theta = \theta_0\dd z_j^3$. Since $\Theta$ is meromorphic at $p_j$, $\theta_0 = \mathrm{const}\cdot z_j^{n_j}$ with $n_j>-\infty$. A. Haydys \cite{isolated_SKS-Haydys2015} showed that there holds a PDE for an SKS metric
    \begin{equation*}
        \frac{\Delta_{\mathbb{R}^2} U}{e^{2U}} = 16\left|\theta_0\right|^2
    \end{equation*}
    Hence, locally on $U_j\setminus\{p_j\}$,
    \begin{equation*}
        K\dd A = \frac{1}{2}e^U\Delta U \dd A
        = 8 e^{3U}\left|\theta_0\right|^2\dd A 
        = 8 e^{2U}\left|\theta_0\right|^2 \left(\frac{\mathrm{i}}{2}\dd z_j \wedge \dd \Bar{z}_j\right).
    \end{equation*}
    Then from the estimation in Lemma \ref{A._local_model_SKS}, around each $p_j$, $K\dd A$ is asymptotical to $\frac{\dd r\wedge\dd \theta}{r(\log r)^2}$ or $\frac{\dd r \wedge \dd\theta}{r^{2(\beta_j - n_j)-1}}$ with $2(\beta_j - n_j)-1 <1$. 
    Thus $K$ is $L^1$ integrable w.r.t. $\dd A$.
\qquad \qquad\qquad\qquad\qquad\qquad\qquad\qquad \qquad \qquad\qquad\qquad\qquad\qquad\qquad\qquad \qquad \qquad \qquad
\qquad \qquad
\boxed{}

\bibliographystyle{plain}
\bibliography{ref}

\begin{thebibliography}{1}

\bibitem{Atiyah-Lebrun2013}
MICHAEL ATIYAH and CLAUDE LEBRUN.
\newblock Curvature, cones and characteristic numbers.
\newblock {\em Mathematical Proceedings of the Cambridge Philosophical Society}, 155(1):13–37, 2013.

\bibitem{Chern_intrinsic-proof-GaussBonnet1944}
Shiing-shen Chern.
\newblock A simple intrinsic proof of the {G}auss-{B}onnet formula for closed {R}iemannian manifolds.
\newblock {\em Ann. of Math. (2)}, 45:747--752, 1944.

\bibitem{Gauss-Bonnet_FXY2024}
Han-bing Fang, Bin Xu, and Bai-rui Yang.
\newblock The {G}auss-{B}onnet formula of a conical metric on a compact {R}iemann surface.
\newblock {\em Chinese Quart. J. Math.}, 39(2):180--184, 2024.

\bibitem{spK-manifolds_Freed1900}
Daniel~S. Freed.
\newblock Special {K}\"ahler manifolds.
\newblock {\em Comm. Math. Phys.}, 203(1):31--52, 1999.

\bibitem{isolated_SKS-Haydys2015}
A.~Haydys.
\newblock Isolated singularities of affine special {K}{\"a}hler metrics in two dimensions.
\newblock {\em Comm. Math. Phys.}, 340(3):1231--1237, 2015.

\bibitem{Haydys-Xu_sks2020}
Andriy Haydys and Bin Xu.
\newblock Special {K}\"ahler structures, cubic differentials and hyperbolic metrics.
\newblock {\em Selecta Math. (N.S.)}, 26(3):Paper No. 37, 21, 2020.

\bibitem{Hormander_PDE-I-1983}
Lars H\"ormander.
\newblock {\em The analysis of linear partial differential operators. {I}}, volume 256 of {\em Grundlehren der mathematischen Wissenschaften [Fundamental Principles of Mathematical Sciences]}.
\newblock Springer-Verlag, Berlin, 1983.
\newblock Distribution theory and Fourier analysis.

\bibitem{Mcmullen_GaussBonnet-cone-manifolds2017}
Curtis~T. McMullen.
\newblock The {G}auss-{B}onnet theorem for cone manifolds and volumes of moduli spaces.
\newblock {\em Amer. J. Math.}, 139(1):261--291, 2017.

\bibitem{Troyanov_compact-surface-conical-singularities1991}
Marc Troyanov.
\newblock Prescribing curvature on compact surfaces with conical singularities.
\newblock {\em Trans. Amer. Math. Soc.}, 324(2):793--821, 1991.

\end{thebibliography}

\end{document}